\documentclass{article}

 \usepackage[sglblindworkshop]{neurips_2025}

\workshoptitle{}

\usepackage[utf8]{inputenc} 
\usepackage[T1]{fontenc}    
\usepackage{amsmath}        
\usepackage{amsthm}         
\usepackage{thmtools,thm-restate} 
\usepackage{url}            
\usepackage{booktabs}       
\usepackage{amsfonts}       
\usepackage{nicefrac}       
\usepackage{microtype}      
\usepackage{xcolor}         
\usepackage{graphicx}       
\usepackage{amssymb}
\usepackage{esint}
\usepackage{geometry}
\usepackage{fancyhdr}
\usepackage{mathrsfs}
\usepackage{mathtools}
\usepackage{tikz}           
\usepackage{hyperref}       

\newcommand{\vs}{\vspace{0.2cm}}
\newcommand{\no}{\noindent}
\newcommand{\nLine}{\vs \no}

\newcommand{\C}[1]{\mathcal{#1}}

\newcommand{\lr}[1]{\left({#1}\right)}
\newcommand{\eps}{\varepsilon}

\title{Prioritizing Recurrent Services}

\author{%
  Lin (Franklin) Feng \\
  Stanford University\\
  \And
  Yue Hu\\
  Stanford University\\
  \And
  Xu Kuang \\
  Stanford University\\
}

\begin{document}
\nolinenumbers

\maketitle

\begin{abstract}
We study optimal scheduling in multi-class queueing systems with reentrance, where jobs may return for additional service after completion. Such reentrance creates feedback loops that fundamentally alter congestion dynamics and challenge classical scheduling results. We model two distinct dimensions of the reentrance behavior, the probability of return and the speed of return, and show that their product, the {\it effective return rate}, is the key statistic that governs optimal priorities. Our main result establishes a dichotomy: when the effective return rate of the smaller job class (the class with lower expected total workload) is lower, a fixed priority rule is optimal; when it is higher, fixed rules are suboptimal and the optimal policy must be state dependent. This characterization clarifies how reentrance changes the externalities that jobs impose on one another and provides structural guidance for designing scheduling policies.
\end{abstract}

\section{Introduction}\label{Sec:Intro}
Service systems across many domains routinely face recurrent demand. In healthcare, patients return for follow-up visits; in call centers, customers make repeated support calls; and in professional services, clients revisit financial advisors. Service delivered today can therefore generate future arrivals, creating feedback loops that fundamentally reshapes system dynamics. Queueing models such as Erlang-R \citep{yom2014erlang} explicitly capture these returns and show that ignoring them can lead to substantial biases in both analysis and decision making.

Our goal is to better control service systems with reentrance through smart prioritization and scheduling. Two observations motivate our study. First, recent advancements in predictive analytics and artificial intelligence make it increasingly feasible to classify jobs by their likelihood and timing of return. For example, in tech support, most setup or installation calls are resolved on the first attempt, but if they require a return, the follow-up typically happens within hours, whereas subscription clients almost always require support again, though only after weeks or months. Second, classical scheduling policies, such as the $c\mu$-rule \citep{mandelbaum2004scheduling} and shortest remaining processing time (SRPT) scheduling \citep{dong2021srpt}, have been shown to be effective because they prioritize jobs that impose the least {\it externality} on others. This raises a natural question: in the presence of reentrance, where jobs differ in both their probability and speed of return, how should externality be quantified, and what scheduling policies remain effective?

We address this question using a fluid queueing model where a system manager dynamically allocates service capacity between two job classes (indexed by $i \in \{1, 2\}$). The classes differ in their return probability ($r_i$) and return rate ($\gamma_i$). The system consists of a primary queue, where jobs line up for service and incur holding costs while waiting, and a virtual return queue, which tracks jobs that will reenter the system in the future. The objective is to minimize total holding costs in the primary queue. Our main research question is how optimal scheduling priorities should be structured in the presence of reentrance. We focus on whether the classical intuition of prioritizing ``smaller'' jobs with lower expected workload remains valid under reentrance, and whether optimal decisions are governed by fixed priority rules or require adaptation to system congestion.

\begin{figure}[ht]
\centering

\begin{tikzpicture}[scale=1.5]
    \def\rA{0.4}
    \def\rB{0.8}
    \begin{scope}
      \clip (0,0) rectangle (3.0,3.0);
      \fill[blue!10] (0,0) -- (3,{(\rA/\rB)*3}) -- (3,0) -- cycle;
      \fill[red!10] (0,0) -- (3,{(\rA/\rB)*3}) -- (3,3) -- (0,3) -- cycle;
    \end{scope}
    
    \draw[thick, ->] (-0.5,0) -- (3.3,0) node[right] {$\gamma_1$};
    \draw[thick, ->] (0,-0.5) -- (0,3.3) node[left] {$\gamma_2$};

    \draw[thick, red!30!black] (0,0) -- ({3}, {(\rA/\rB)*3});

    \node[red!30!black, anchor=west] at (3,1.5) {$\kappa_1=\kappa_2$};
    \node[red!70!black, anchor=west] at (3,2.5) {Fixed priority policy optimal};
    \node[red!70!black, anchor=west] at (3,2.3) {$\kappa_1\leq\kappa_2$};
    \node[blue!70!black, anchor=west] at (3,0.7) {State-dependent policy optimal};
    \node[blue!70!black, anchor=west] at (3,0.5) {$\kappa_1>\kappa_2$};
    
\end{tikzpicture}
\caption{Structure of the optimal scheduling policy}
\label{fig:two_regime}
\end{figure}
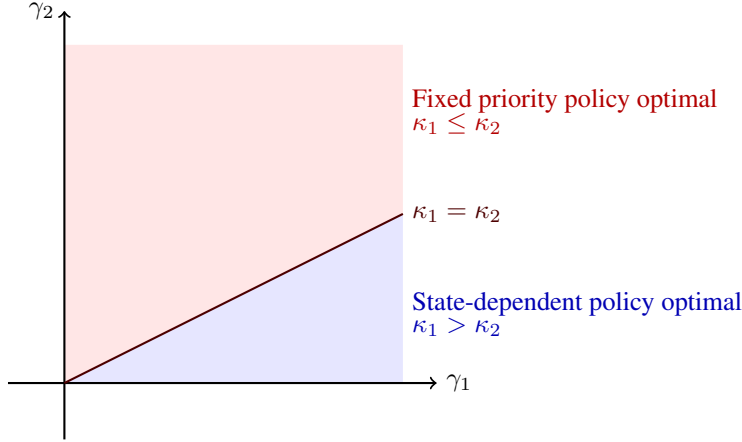

Our main result provides a sharp characterization of the optimal scheduling policy; see Figure \ref{fig:two_regime}.
Assume without loss of generality that $r_1 < r_2$, so that class 1 represents the smaller jobs with lower expected total workload. For each class $i$, define the effective return rate as $\kappa_i = r_i \gamma_i$, which jointly captures how likely and how quickly work returns.
We find that:
\begin{enumerate}
    \item If $\kappa_1 \leq \kappa_2$, a fixed priority policy is optimal: class 1 (the smaller jobs) should always be prioritized, independently of the system’s congestion level.
    \item If $\kappa_1 > \kappa_2$, no fixed priority policy is optimal. In this case, class 1 should be prioritized under heavy congestion, but priority should shift to class 2 when the system is lightly loaded.
\end{enumerate}

A key insight of our analysis is that reentrance changes the way jobs impose externalities on one another. In classical multi-class queues without returns, smaller jobs are always favored because they clear faster and permanently relieve congestion.
With reentrance, however, serving jobs that are likely to return soon can offset this benefit. The effective return rate $\kappa_i = r_i\gamma_i$ captures this feedback and becomes the critical (though not exclusive) measure of class $i$'s externality.
When $\kappa_1 \leq \kappa_2$, the benefit of completing the smaller job still dominates, so always prioritizing class 1 remains optimal.
When $\kappa_1 > \kappa_2$, the optimal policy switches depending on the system load. When the system is already heavily loaded, the overriding concern is to {\em drain work quickly} to reduce holding costs. Class 1 jobs are smaller (since $r_1 < r_2$), serving them clears backlog faster and immediately relieves congestion, even though they are prone to return soon. In this regime, the benefit of faster clearance outweighs the risk of quick reentries. Under light congestion, however, the immediate backlog is less pressing, and the main concern shifts to the {\em future workload generated by today’s service}. 
It is therefore better to prioritize class 2 with a lower effective return rate.
In short, the effective return rate captures long-run feedback externalities, while the optimal policy balances them ($\kappa_1$ vs. $\kappa_2$) against the immediate externalities of leaving work unfinished ($r_1$ vs. $r_2$), which explains why fixed priority rules suffice in some regimes but fail in others.

\section{Related Literature}

First, our work is related to the literature on optimal scheduling in multi-class queues. One rich line of work emphasizes the power of simple index rules, such as the $c\mu$-rule and its extensions, and shows that these static priority rules can yield (near-)optimal performance in a wide range of queueing models \citep{smith1956various, van1995dynamic, mandelbaum2004scheduling, atar2010cmu, long2020dynamic}. 
Complementing index policies, state-dependent rules that favor short jobs, such as SRPT and service-age-based variants, are known to minimize delay in a range of models \citep{schrage1966queue, scully2018soap, dong2021srpt, ibrahim2026shortest}. We contribute to this literature by showing how reentrance changes the notion of externality that underpins these rules. In our model, the product of return probability and return speed becomes the key statistic that governs whether fixed priority rules suffice or whether state dependence is essential.

Second, our work is related to the stream of literature on reentrant service systems. Classical Erlang-R models demonstrate that ignoring reentrance can lead to biased performance estimates and miscalibrated staffing, especially in healthcare \citep{yom2014erlang, armony2015patient}. More recent work shows that heterogeneous return dynamics can alter system stability and equilibrium \citep{barjesteh2021multiclass}. In addition, reentrance has been incorporated into diverse applications, including emergency department staffing with time-varying physician productivity \citep{ouyang2021emergency}, post-discharge hospital readmission prevention \citep{chan2025dynamic}, community corrections placement \citep{gao2025stopping}, and customer-agent interactions in contact centers \citep{daw2025co}. Collectively, these studies model reentrance as an important operational feature, but relatively little is known about how to optimally schedule recurrent jobs. We address this gap directly by showing how reentrance reshapes the notion of externalities in scheduling theory and by characterizing the structure of the optimal scheduling policy.

\section{The Model}\label{sec:model}
We consider a two-class fluid queueing system with reentrant jobs. The system is closed, with no external arrivals, and the objective is to optimally clear all existing workload.
For each class $i$, work first enters a {\em primary queue}, where it receives service. After service completion, a fraction $r_i \in (0, 1)$ of the work reenters the system by joining a virtual {\em return queue}. Jobs in the return queue depart at rate $\gamma_i > 0$, at which point they rejoin the primary queue.
Without loss of generality, we assume $r_1 < r_2$, so that class 1 represents ``smaller'' jobs with a lower expected total work once the reentrance probability is taken into account. 

The system manager has a total service capacity $\mu>0$, and dynamically allocates it between the two classes. Specifically, the system manager determines allocations $u(t)=(u_1(t), u_2(t))$, subject to $u_i(t)\geq 0$ and $u_1(t) + u_2(t)\leq 1$, $i \in \{1, 2\}$, $t \geq 0$. 
If class $i$ receives a fraction $u_i(t)$ of capacity at time $t$, its primary queue is depleted at rate $\mu u_i(t)$. 

At time $t$, let $q_i^p(t)$ denote the primary queue length of class $i$, $q_i^r(t)$ the amount of work in the return queue that will reenter in the future, and $q(t)=(q_1^p(t),q_1^r(t), q_2^p(t),q_2^r(t))$ the full system state. The system dynamics are given by
\begin{equation}
    \dot q_i^p(t) = -\mu u_i(t) + \gamma_i q_i^r(t), \quad \dot q_i^r(t) = r_i \mu u_i(t) - \gamma_i q_i^r(t), \quad i = 1, 2.\label{eq:dynamics}
\end{equation}

We call a mapping $\psi$ an \emph{admissible policy} if it prescribes allocations $u(t)=\psi(q(t))$ that maintain all queues nonnegative. A sample path of allocations $u$ is an \emph{optimal trajectory} for given parameters and an initial condition if, among all admissible allocations, it achieves the minimal cumulative holding cost incurred by primary queues over a sufficiently long horizon $T$:
\begin{equation*}
\int_0^T \left( q_1^p(t) + q_2^p(t) \right) \,dt.
\end{equation*}
A policy $\psi$ is said to be \emph{optimal} for a given parameters if the resulting allocation trajectory $u(t)=\psi(q(t))$, $t\geq 0$, is an optimal trajectory for any initial condition. In general, we are interested in understanding how to construct optimal policies. 

We say that the server \emph{prioritizes} class $i\in\{1,2\}$ at state $q$ if  
(i) when $q_i^p>0$, the server devotes full capacity to class $i$;  
and (ii) when $q_i^p=0$, the server allocates enough capacity to keep $q_i^p$ empty.  

A policy $\psi$ is a \emph{fixed priority policy} if the server prioritizes one particular class at all states $q$, and a \emph{state-dependent policy} otherwise. To emphasize, by state-dependent policy we mean one that is strictly state-dependent; fixed priority rules are not included as a special case.

\section{Main Results}
We now establish how reentrance shapes the structure of the optimal scheduling policy. The key determinant is the {\em effective return rate} $\kappa_i = r_i\gamma_i$, which reflects not only how much future workload class $i$ generates after service, but also the speed at which this workload returns to the system. Comparing $\kappa_1$ and $\kappa_2$ yields two distinct regimes.

\begin{restatable}{theorem}{OptPolicyThm}\label{thm:opt_policy}
The optimal scheduling policy satisfies: 
\begin{enumerate}
 \item If $\kappa_1 \le \kappa_2$, a fixed priority policy that always prioritizes class $1$ is optimal.
\item If $\kappa_1 > \kappa_2$, no fixed priority policy is optimal, and the optimal policy is state-dependent.

\end{enumerate}   
\end{restatable}

\begin{figure}[ht]
\centering
\includegraphics[scale = 0.46]{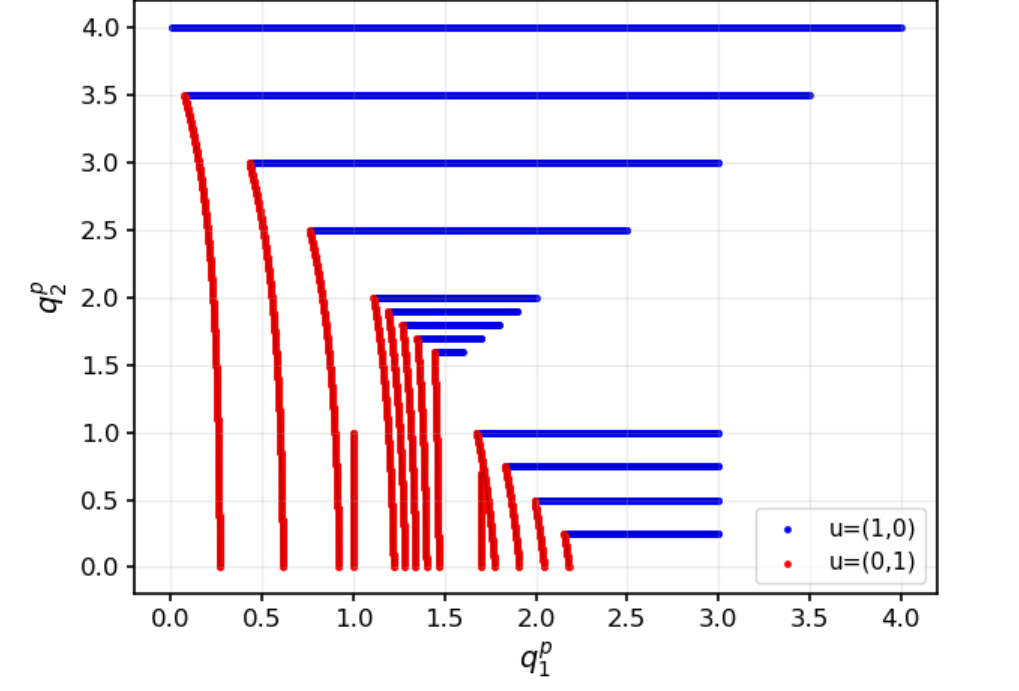}
\vspace{-.2cm}
\caption{{\small Optimal state trajectories when $\kappa_1>\kappa_2$ ($r_1 =0.2$, $r_2 = 0.8$, $\gamma_1=2$, $\gamma_2 = 0.2$, $\mu = 2$, $q^r(0)=0$)}}\label{fig:traj}
\end{figure}

We complement Theorem~\ref{thm:opt_policy} with the following important observation, supported by extensive numerical experiments.
In the regime where $\kappa_1 > \kappa_2$, the optimal state-dependent policy exhibits at most one switch along any trajectory. Specifically, the policy prioritizes class 1 under heavy congestion and transitions to class 2 as the system clears.
We next present a set of numerical experiments with parameters from the region $\kappa_1 > \kappa_2$ where no fixed priority policy is optimal. Figure \ref{fig:traj} plots the evolution of the state $(q_1^p,q_2^p)$ under the optimal policy, starting from different initial conditions. We make two consistent observations: (i) when both queues are heavily loaded, the policy prioritizes class 1, but as the system clears, priority shifts to class 2; and (ii) along every trajectory, at most one switch of priority occurs.
In other words, smaller jobs dominate prioritization under high congestion, but their rapid returns eventually make it more efficient to prioritize larger jobs as the system clears.

Importantly, our results illustrate how reentrance changes the nature of externalities in scheduling. In classical queues without reentrance, smaller jobs always impose less externality, since they leave the system quickly and permanently free capacity. With reentrance, completing a job may regenerate demand, so a higher effective return rate $\kappa_i$ reflects a stronger {\em future} externality. 
However, the optimal rule is not determined by comparing $\kappa_1$ and $\kappa_2$ alone. Under heavy congestion, the main externality is the {\em immediate} delay from large backlogs, and prioritizing the smaller jobs (class 1) best alleviates this burden, even if they are more likely to return soon. 
Under light congestion, the immediate backlog is less pressing, and the dominant externality becomes the feedback created by future returns. In that regime, it is optimal to prioritize the class with the smaller $\kappa_i$.
In a nutshell, $\kappa_i$ quantifies future externalities, but the optimal prioritization balances them against the more myopic externality of delaying work clearance. This tradeoff explains why fixed priority rules fail when $\kappa_1 > \kappa_2$ and state dependence becomes necessary.

To conclude, we numerically demonstrate the value of state-dependent policies by comparing the objective values of the optimal state-dependent policy with those of the two fixed priority rules. In Table \ref{tab:gaps}, the ``FP-1 gap'' column reports the relative performance loss from always prioritizing class 1, while the ``FP-2 gap'' column does the same for class 2. The results show that giving fixed priority to class 1 can be close to optimal when the parameters are near the boundary where $\kappa_1=\kappa_2$, but becomes increasingly suboptimal as class 2 returns more slowly. Conversely, fixed priority to class 2 performs poorly near the boundary but improves steadily as its return rate decreases.

\begin{table}[ht]
\centering
\renewcommand{\arraystretch}{1.1}
\begin{tabular}{@{\hskip 0.4cm}c@{\hskip 0.4cm}|@{\hskip 0.4cm}c@{\hskip 0.4cm}|@{\hskip 0.4cm}c@{\hskip 0.4cm}}
\hline
$\gamma_2$ & FP-1 gap & FP-2 gap \\
\hline
0.05 & 7.35\% & 0.01\% \\
0.10 & 5.02\% & 0.01\% \\
0.20 & 0.89\% & 0.67\% \\
0.30 & 0.03\% & 4.68\% \\
0.40 & 0.01\% & 10.16\% \\
0.50 & 0.00\% & 16.01\% \\
\hline
\end{tabular}
\vspace{.2cm}
\caption{{\small Improvement over fixed priority rules ($r_1 =0.2$, $r_2 = 0.8$, $\gamma_1=2$, $\mu = 2$, $q^p(0)=2$, $q^r(0)=0$)}}
\label{tab:gaps}
\end{table}


\newpage
\bibliography{ref}

\begin{thebibliography}{16}
\providecommand{\natexlab}[1]{#1}
\providecommand{\url}[1]{\texttt{#1}}
\expandafter\ifx\csname urlstyle\endcsname\relax
  \providecommand{\doi}[1]{doi: #1}\else
  \providecommand{\doi}{doi: \begingroup \urlstyle{rm}\Url}\fi

\bibitem[Armony et~al.(2015)Armony, Israelit, Mandelbaum, Marmor, Tseytlin, and Yom-Tov]{armony2015patient}
Mor Armony, Shlomo Israelit, Avishai Mandelbaum, Yariv~N Marmor, Yulia Tseytlin, and Galit~B Yom-Tov.
\newblock On patient flow in hospitals: A data-based queueing-science perspective.
\newblock \emph{Stochastic systems}, 5\penalty0 (1):\penalty0 146--194, 2015.

\bibitem[Atar et~al.(2010)Atar, Giat, and Shimkin]{atar2010cmu}
Rami Atar, Chanit Giat, and Nahum Shimkin.
\newblock The c$\mu$/$\theta$ rule for many-server queues with abandonment.
\newblock \emph{Operations Research}, 58\penalty0 (5):\penalty0 1427--1439, 2010.

\bibitem[Barjesteh and Abouee-Mehrizi(2021)]{barjesteh2021multiclass}
Nasser Barjesteh and Hossein Abouee-Mehrizi.
\newblock Multiclass state-dependent service systems with returns.
\newblock \emph{Naval Research Logistics (NRL)}, 68\penalty0 (5):\penalty0 631--662, 2021.

\bibitem[Chan et~al.(2025)Chan, Huang, and Sarhangian]{chan2025dynamic}
Timothy~CY Chan, Simon~Y Huang, and Vahid Sarhangian.
\newblock Dynamic control of service systems with returns: Application to design of postdischarge hospital readmission prevention programs.
\newblock \emph{Operations Research}, 73\penalty0 (4):\penalty0 2242--2263, 2025.

\bibitem[Daw et~al.(2025)Daw, Castellanos, Yom-Tov, Pender, and Gruendlinger]{daw2025co}
Andrew Daw, Antonio Castellanos, Galit~B Yom-Tov, Jamol Pender, and Leor Gruendlinger.
\newblock The co-production of service: Modeling services in contact centers using hawkes processes.
\newblock \emph{Management Science}, 71\penalty0 (3):\penalty0 2635--2656, 2025.

\bibitem[Dong and Ibrahim(2021)]{dong2021srpt}
Jing Dong and Rouba Ibrahim.
\newblock Srpt scheduling discipline in many-server queues with impatient customers.
\newblock \emph{Management Science}, 67\penalty0 (12):\penalty0 7708--7718, 2021.

\bibitem[Gao et~al.(2025)Gao, Shi, and Kong]{gao2025stopping}
Xiaoquan Gao, Pengyi Shi, and Nan Kong.
\newblock Stopping the revolving door: Mdp-based decision support for community corrections placement.
\newblock \emph{Available at SSRN 4672337}, 2025.

\bibitem[Ibrahim and Dong(2026)]{ibrahim2026shortest}
Rouba Ibrahim and Jing Dong.
\newblock Shortest-job-first scheduling in many-server queues with impatient customers and noisy service-time estimates.
\newblock \emph{Operations Research}, 2026.

\bibitem[Long et~al.(2020)Long, Shimkin, Zhang, and Zhang]{long2020dynamic}
Zhenghua Long, Nahum Shimkin, Hailun Zhang, and Jiheng Zhang.
\newblock Dynamic scheduling of multiclass many-server queues with abandonment: The generalized c$\mu$/h rule.
\newblock \emph{Operations Research}, 68\penalty0 (4):\penalty0 1218--1230, 2020.

\bibitem[Mandelbaum and Stolyar(2004)]{mandelbaum2004scheduling}
Avishai Mandelbaum and Alexander~L Stolyar.
\newblock Scheduling flexible servers with convex delay costs: Heavy-traffic optimality of the generalized c$\mu$-rule.
\newblock \emph{Operations Research}, 52\penalty0 (6):\penalty0 836--855, 2004.

\bibitem[Ouyang et~al.(2021)Ouyang, Liu, and Sun]{ouyang2021emergency}
Huiyin Ouyang, Ran Liu, and Zhankun Sun.
\newblock Emergency department modeling and staffing: Time-varying physician productivity.
\newblock \emph{Available at SSRN 3963226}, 2021.

\bibitem[Schrage and Miller(1966)]{schrage1966queue}
Linus~E Schrage and Louis~W Miller.
\newblock The queue m/g/1 with the shortest remaining processing time discipline.
\newblock \emph{Operations Research}, 14\penalty0 (4):\penalty0 670--684, 1966.

\bibitem[Scully et~al.(2018)Scully, Harchol-Balter, and Scheller-Wolf]{scully2018soap}
Ziv Scully, Mor Harchol-Balter, and Alan Scheller-Wolf.
\newblock Soap: One clean analysis of all age-based scheduling policies.
\newblock \emph{Proceedings of the ACM on Measurement and Analysis of Computing Systems}, 2\penalty0 (1):\penalty0 1--30, 2018.

\bibitem[Smith et~al.(1956)]{smith1956various}
Wayne~E Smith et~al.
\newblock Various optimizers for single-stage production.
\newblock \emph{Naval Research Logistics Quarterly}, 3\penalty0 (1-2):\penalty0 59--66, 1956.

\bibitem[Van~Mieghem(1995)]{van1995dynamic}
Jan~A Van~Mieghem.
\newblock Dynamic scheduling with convex delay costs: The generalized c| mu rule.
\newblock \emph{The Annals of Applied Probability}, pages 809--833, 1995.

\bibitem[Yom-Tov and Mandelbaum(2014)]{yom2014erlang}
Galit~B Yom-Tov and Avishai Mandelbaum.
\newblock Erlang-r: A time-varying queue with reentrant customers, in support of healthcare staffing.
\newblock \emph{Manufacturing \& Service Operations Management}, 16\penalty0 (2):\penalty0 283--299, 2014.

\end{thebibliography}


\newpage
\appendix

\section{Proof of Theorem \ref{thm:opt_policy}}
In this section, we present the proof of Theorem \ref{thm:opt_policy}. Section \ref{subsec:proof_opt_policy_FP} shows that the fixed-priority policy that always prioritizes class 1 is optimal when $\kappa_1 \leq \kappa_2$. Section \ref{subsec:proof_opt_policy_SDP} establishes that no fixed-priority policy can be uniformly optimal for all initial states $q_0$ when $\kappa_1 > \kappa_2$, using proof by contradiction.

\subsection{Proof of Theorem \ref{thm:opt_policy} Part (1)}\label{subsec:proof_opt_policy_FP}
To prove Part (1) of Theorem \ref{thm:opt_policy}, we proceed in two steps. First, we rewrite the control problem by expressing the state trajectories $q_i^p(t)$ and $q_i^r(t)$ in integral form. This representation yields a convenient kernel formulation of the objective, highlighting how the control enters linearly over time. Second, we apply Pontryagin’s Minimum Principle. By analyzing the Hamiltonian and its coefficients, we show that whenever $\kappa_1 \leq \kappa_2$, the Hamiltonian is minimized by allocating as much capacity as feasibly possible to class 1 at every state. This corresponds exactly to the fixed priority policy that always prioritizes class 1, thereby proving its optimality.

We now begin the proof by expressing $q_i^p(t)$ and $q_i^r(t)$ in integral form. By \eqref{eq:dynamics}, for any feasible allocation $u=(u_1,u_2)$ (i.e., allocation that satisfies $0\leq u_i\leq 1$ and $u_1+u_2\leq 1$), we obtain
\begin{align*}
    q_i^r(t)&=e^{-\gamma_i t}q_{i,0}^r+r_i\mu\int_0^t e^{-\gamma_i(t-s)}u_i(s)\,ds,\\
    q_i^p(t)&=q_{i,0}^p+\int_0^t -\mu u_i(s)+\gamma_iq_i^r(s)\,ds\\
    &=q_{i,0}^p-\mu \int_0^tu_i(s)\,ds+\int_0^t\gamma_iq_i^r(s)\,ds\\
    &=q_{i,0}^p-\mu \int_0^tu_i(s)\,ds+(1-e^{-\gamma_i t})q_{i,0}^r+\gamma_i r_i\mu \int_0^t\left[\int_s^t e^{-\gamma_i(x-s)}\,dx\right]u_i(s)\,ds\\
    &=q_{i,0}^p-\mu \int_0^t u_i(s)\,ds+(1-e^{-\gamma_i t})q_{i,0}^r+\gamma_i r_i\mu \int_0^t \frac{1}{\gamma_i}\left[1-e^{-\gamma_i(t-s)}\right]u_i(s)\,ds.
\end{align*}
Define the kernel $k_i(\tau)=\mu[(r_i-1)-r_i e^{-\gamma_i \tau}]$, $\tau\geq 0$, so that the primary queue length admits the compact representation
\begin{align*}
    q_i^p(t)=q_{i,0}^p+(1-e^{-\gamma_i t})q_{i,0}^r+\int_0^t k_i(t-s)u_i(s)\,ds.
\end{align*}
Next, let
\begin{align*}
    K_i(\tau)&=\int_0^{\tau}k_i(x)\,dx\\
    &=\mu\left[(r_i-1)\tau-\frac{r_i}{\gamma_i}(1-e^{-\gamma_i\tau})\right]
\end{align*}
denote the cumulative kernel. Using Fubini's Theorem, the finite-horizon objective can be expressed as
\begin{align*}
    \int_0^T q_1^p(t)+q_2^p(t)\,dt&=\sum_{i=1}^2\int_0^T\left[q_{i,0}^p+(1-e^{-\gamma_i t})q_{i,0}^r\right]\,dt+\sum_{i=1}^2\int_0^T\int_0^t k_i(t-s)u_i(s)\,ds\,dt\\
    &=C_T(q_0)+\sum_{i=1}^2\int_0^T u_i(s)\left[\int_{s}^T k_i(t-s)\,dt\right]ds\\
    &=C_T(q_0)+\sum_{i=1}^2\int_0^Tu_i(s)K_i(T-s)\,ds,
\end{align*}
where $C_T(q_0)$ depends only on the initial state $q_0$ and horizon $T$. Therefore, minimizing the cumulative holding cost is equivalent to solving
\begin{align*}
    \min_{u(\cdot)}\sum_{i=1}^2\int_0^Tu_i(s)K_i(T-s)\,ds,
\end{align*}
which reveals the problem as a linear functional of the control trajectory with weights $K_i(T-s)$. Importantly, the weight function $K_i(\cdot)$ has the following properties,
\begin{restatable}{lemma}{KiPropsLem}\label{lem:K_i_props}
    Assume that $0<r_1<r_2<1$, for all $\tau\geq 0$, $K_i(\tau)\leq 0$. In particular, when $\kappa_1\leq \kappa_2$, $K_1(\tau)\leq K_2(\tau)$.
\end{restatable}
\textit{Proof. }See Section \ref{subsec:proof_supp_lem}.

Recall from Section~\ref{sec:model} that a policy $\psi$ is admissible if the allocations $u(t)=\psi(q(t))$ it prescribes keep all queues nonnegative. In particular, when the primary queue of class $i$ is empty ($q_i^p(t)=0$), we require $\dot{q}_i^p(t)\geq 0$, which implies the constraint $\mu u_i(t)\leq \gamma_i q_i^r(t)$. Collecting these conditions, the admissible allocation set at state $q$ is defined as
\begin{align}
    \C{U}(q):=\left\{u: u_1,u_2\geq 0, u_1+u_2\leq 1, \text{if }q_i^p=0, \text{then }\mu u_i\leq\gamma_iq_i^r\text{ for }i=1,2\right\}.\label{eq:CU_q}
\end{align}
We now apply Pontryagin's Principle to show that the fixed-priority policy prioritizing class 1 is optimal. For a fixed horizon $T>0$, the Hamiltonian is
\begin{align*}
    \C{H}(q,u,\lambda)=\sum_{i=1}^2q_i^p+\sum_{i=1}^2\lambda_i^p(-\mu u_i+\gamma_i q_i^r)+\sum_{i=1}^2\lambda_i^r(r_i\mu u_i-\gamma_iq_i^r).
\end{align*}
The costates satisfy $\dot{\lambda}_i^p(t)=-1$, $\dot{\lambda}_i^r(t)=-\gamma_i(\lambda_i^p(t)-\lambda_i^r(t))$, $\lambda_i^p(T)=\lambda_i^r(T)=0$. Solving backward gives us
\begin{align*}
    \lambda_i^p(t)=T-t,\quad \lambda_i^r(t)=T-t-\frac{1}{\gamma_i}(1-e^{-\gamma_i(T-t)}).
\end{align*}
So the coefficient of $u_i$ in $\C{H}$ is
\begin{align*}
    \mu(r_i\lambda_i^r(t)-\lambda_i^p(t))=\mu\left[ r_i(T-t)-\frac{r_i}{\gamma_i}(1-e^{-\gamma_i(T-t)}) -(T-t)\right]=K_i(T-t).
\end{align*}
From Lemma \ref{lem:K_i_props}, for every $t\in[0,T)$, $K_1(T-t)\leq K_2(T-t)\leq 0$. By Pontryagin's Principle, an optimal control $u^*(\cdot)$ must minimize the Hamiltonian pointwise over $\C{U}(q^*(t))$, the admissible control set given in \eqref{eq:CU_q} to ensure control feasibility and state nonnegativity. Thus, the optimal allocation is
\begin{align*}
    u_1^*(t)=\begin{cases}
        1,&\text{if }q_1^p(t)>0,\\
        \min\{1,\gamma_1q_1^r(t)/\mu\}&\text{if }q_1^p(t)=0.
    \end{cases}\quad
    u_2^*(t)=\begin{cases}
        1-u_1^*(t),&\text{if }q_2^p(t)>0,\\
        \min\{1-u_1^*(t),\gamma_2q_2^r(t)/\mu\}&\text{if }q_2^p(t)=0.
    \end{cases}
\end{align*}
That is, the optimal policy allocates full capacity to class 1 whenever it has positive workload, and allocate just enough capacity to it to keep $q_1$ at 0 while maintaining feasibility. On the other hand, any leftover capacity will be given to class 2 whenever it has positive workload. This corresponds precisely to the fixed priority policy introduced in Section \ref{sec:model}. Therefore, when $\kappa_1\leq \kappa_2$, the fixed priority policy that always prioritizes class 1 is optimal and obtains minimal cost under any $T$.\qed

\subsection{Proof of Theorem \ref{thm:opt_policy} Part (2)}\label{subsec:proof_opt_policy_SDP}
Part (2) of Theorem \ref{thm:opt_policy} requires showing that no fixed-priority policy is optimal when $\kappa_1>\kappa_2$. Our strategy is to demonstrate that the two fixed-priority policies reverse their cost ranking under different initial conditions. Consider the family of initial states $q_0^{(\eps)}=(q_1^p,q_2^p,q_1^r,q_2^r)=(\eps,\eps,0,0)$ with $\eps>0$. For any parameters $(r_1,r_2,\gamma_1,\gamma_2)$ such that $r_1<r_2$ and $\kappa_1>\kappa_2$, we show that as $\eps\to 0^+$, the policy that always prioritizes class 2 achieves a strictly lower cost than the policy that always prioritizes class 1. In contrast, as $\eps\to\infty$, the inequality reverses. This proves that no single fixed-priority policy can minimize cost for all initial states $q_0$ in this regime.

Before analyzing the inequalities in the two asymptotic regimes, we first derive several key quantities needed for the proof. Suppose we implement the policy that prioritizes class $i$ exclusively. Let $t_i(\eps)$ denote the first time at which the initial backlog of class $i$ is cleared (that is, when $q_i^p(t)$ reaches zero). Then, on interval $[0,t_i(\eps))$, the primary and return queue lengths can be expressed in integral form as follows
\begin{align*}
    q_i^r(t)&=r_i\mu \int_0^te^{-\gamma_i(t-s)}\,ds=\frac{r_i\mu }{\gamma_i}(1-e^{-\gamma_i t}),\\
    q_i^p(t)&=\eps+\mu \int_0^t\left[(r_i-1)-r_ie^{-\gamma_i(t-s)}\right]\,ds=\eps-\mu \left[(1-r_i)t+\frac{r_i}{\gamma_i}(1-e^{-\gamma_i t})\right].
\end{align*}
Therefore, $t_i(\eps)$ is the unique solution to
\begin{align}
    \eps=\mu \left[ (1-r_i)t_i(\eps)+\frac{r_i}{\gamma_i}(1-e^{-\gamma_i t_i(\eps)}) \right].\label{eq:ti_eps}
\end{align}
Let $I_i(\eps)$ denote the total cost accumulated by class $i$ over the interval $[0,t_i(\eps))$,
\begin{align*}
    I_i(\eps)=\int_0^{t_i(\eps)}q_i^p(t)\,dt
    =\eps t_i-\frac{(1-r_i)\mu }{2} t_i^2(\eps)-\frac{r_i\mu }{\gamma_i}\lr{t_i-\frac{1-e^{-\gamma_i t_i(\eps)}}{\gamma_i}}.
\end{align*}
Once the backlog of class $i$ is cleared, the server continues to allocate a fraction $\gamma_i q_i^r(t)/\mu$ of its capacity to class $i$ to keep its primary queue at zero, while the remaining capacity is devoted to the class $j$, which we denoted as the class that is not prioritized.

We now turn to the behavior of class $j$ once the backlog of class $i$ has been cleared. Define $\bar{t}_j(\eps)$ as the first time when the backlog of class $j$ is depleted (that is, when $q_j^p(t)$ reaches zero). For clarity, we divide the trajectory into two stages: 1) \textbf{\emph{Stage A}} as the interval $[0,t_i(\eps))$ during which the server works exclusively on class $i$; 2) \textbf{\emph{Stage B}} as the interval $[t_i(\eps),\bar{t}_j(\eps))$ during which class $i$'s backlog remains at $0$, and all residual capacity is used to serve $j$.

For convenience, define
\begin{align*}
    \lambda_i:=(1-r_i)\gamma_i>0,\quad a_i(\eps):=u_i(t_i)=\frac{\gamma_iq_i^r(t_i)}{\mu},\quad \Delta_j(\eps)=\bar{t}_j(\eps)-t_i(\eps),
\end{align*}
and we write $t_i,\bar{t}_j,a_i,\Delta_j$ in place of $t_i(\eps),\bar{t}_j(\eps),a_i(\eps),\Delta_j(\eps)$, respectively. Let $C_i^{(A)}$ and $C_i^{(B)}$ denote the costs accumulated in stages A and B under the fixed-priority policy that prioritizes class $i$.

From the previous expressions, $C_i^{(A)}=I_i(\eps)+\eps t_i(\eps)$. To compute $C_i^{(B)}$, note that for $t\geq t_i$, the service allocations take the form
\begin{align*}
    u_i(t)=a_ie^{-\lambda_i(t-t_i)},\quad u_j(t)=1-u_i(t).
\end{align*}
At time $t_i$, we have $q_j^r(t_i)=0$ and $q_j^p(t_i)=\eps$. Applying the dynamics \eqref{eq:dynamics}, the return and primary queues of class $j$ evolve as
\begin{align}
    q_j^r(t)&=r_j\mu \int_{t_i}^t e^{-\gamma_j(t-s)}u_j(s)\,ds\nonumber\\
    &=r_j\mu \int_{t_i}^te^{-\gamma_j(t-s)}(1-a_ie^{-\lambda_i(s-t_i)})\,ds\nonumber\\
    &=\frac{r_j\mu }{\gamma_j}\lr{1-e^{-\gamma_j(t-t_i)}}-\frac{r_j\mu a_i}{\gamma_j-\lambda_i}\lr{e^{-\lambda_i(t-t_i)}-e^{-\gamma_j(t-t_i)}},\\
    q_j^p(t)&=\eps-(1-r_j)\mu (t-t_i)+\frac{(1-r_j)\mu a_i}{\lambda_i}\lr{1-e^{-\lambda_i(t-t_i)}}\nonumber\\
    &-\frac{r_j\mu }{\gamma_j}\lr{1-e^{-\gamma_j(t-t_i)}}+\frac{r_j\mu a_i}{\gamma_j-\lambda_i}\lr{e^{-\lambda_i(t-t_i)}-e^{-\gamma_j(t-t_i)}}.\label{eq:qjs_stage_B}
\end{align}
The depletion time $\bar{t}_j(\eps)$ then uniquely solves $q_j^p(t)=0$ in \eqref{eq:qjs_stage_B}. During stage B, only class $j$'s backlog contributes to the cost, so
\begin{align}
    C_i^{(B)}&=\int_{t_i}^{\bar{t}_j}q_j^p(t)\,dt\nonumber\\
    &=\eps\Delta_j-\frac{1}{2}(1-r_j)\mu \Delta_j^2+\frac{(1-r_j)\mu a_i}{\lambda_i}\lr{\Delta_j+\frac{e^{-\lambda_i\Delta_j}-1}{\lambda_i}}\nonumber\\
    &-\frac{r_j\mu }{\gamma_j}\lr{\Delta_j+\frac{e^{-\gamma_j\Delta_j}-1}{\gamma_j}}+\frac{r_j\mu a_i}{\gamma_j-\lambda_i}\lr{\frac{1-e^{-\lambda_i\Delta_j}}{\lambda_i}-\frac{1-e^{-\gamma_j\Delta_j}}{\gamma_j}}.\label{eq:C_iB}
\end{align}
Finally, note that $\dot{q}_j^p(t)+\dot{q}_j^r(t)=(r_j-1)\mu u_j(t)$. Integrating this relation over stage B gives
\begin{align*}
    \int_{t_i}^{\bar{t}_j}u_j(t)\,dt=\frac{\eps-q_j^r(\bar{t}_j)}{(1-r_j)\mu },
\end{align*}
Since during stage B we have $u_i(t)=\gamma_iq_i^r(t)/\mu$ and $u_j(t)=1-u_i(t)$, it follows that
\begin{align}
    \Delta_j-\frac{a_i}{\lambda_i}\lr{1-e^{-\lambda_i\Delta_j}}=\frac{\eps-q_j^r(\bar{t}_j)}{(1-r_j)\mu }.\label{eq:Delta_j}
\end{align}
Of which,
\begin{align}
    q_j^r(\bar{t}_j)=\frac{r_j\mu }{\gamma_j }(1-e^{-\gamma_j \Delta_j})-\frac{r_j\mu a_i}{\gamma_j -\lambda_i}\lr{e^{-\lambda_i\Delta_j}-e^{-\gamma_j \Delta_j}}.\label{eq:qjr_bar_tj}
\end{align}
Equations \eqref{eq:Delta_j} and \eqref{eq:qjr_bar_tj} together provide an implicit characterization of $\Delta_j$, which completes the formulation of $C_i^{(B)}$. The total cost under the fixed-priority policy that prioritizes class $i$ is then given by $C_i^{(A)}+C_i^{(B)}$.

\subsubsection{Small Initial Loads \texorpdfstring{$\eps\to 0^+$}{eps to 0+}}
In this section, we show that when $\eps\to 0^+$, the fixed-priority policy prioritizing class 2 yields a strictly lower cost than the policy that always prioritizes class 1. That is, $C_2^{(A)}+C_2^{(B)}<C_1^{(A)}+C_1^{(B)}$ as $\eps\to 0^+$.

Denote $x=\eps/\mu\to 0^+$. By Taylor's expansion on \eqref{eq:ti_eps} and $a_i(\eps)=\gamma_iq_i^r(t_i)/\mu$, we obtain
\begin{align*}
    t_i&=x+\frac{r_i\gamma_ix^2}{2}+\frac{(\gamma_i)^2r_i(3r_i-1)x^3}{6}+O(x^4),\\
    a_i&=r_i\gamma_ix-\frac{(1-r_i)r_i(\gamma_i)^2x^2}{2}+O(x^3).
\end{align*}
Next, we expand both sides of \eqref{eq:Delta_j} in powers of $x$ up to second order, which requires a corresponding expansion of the terms in $q_j^r(\bar{t}_j)$, for which we use the following lemma.
\begin{restatable}{lemma}{DeltajLem}\label{lem:Delta_j_Ox}
    $\Delta_j=O(x)$.
\end{restatable}
\textit{Proof. }See Section \ref{subsec:proof_supp_lem}.

\nLine Then, consider the expansions
\begin{align*}
    1-e^{-\lambda_i\Delta_j}&=\lambda_i\Delta_j-\frac{\lambda_i^2\Delta_j^2}{2}+O(\Delta_j^3),\\
    1-e^{-\gamma_j\Delta_j}&=\gamma_j\Delta_j-\frac{(\gamma_j)^2}{2}\Delta_j^2+O(\Delta_j^3),\\
    e^{-\lambda_i\Delta_j}-e^{-\gamma_j\Delta_j}&=(1-\lambda_i\Delta_j+\frac{\lambda_i^2}{2}\Delta_j^2)-(1-\gamma_j\Delta_j+\frac{(\gamma_j)^2}{2}\Delta_j^2)+O(\Delta_j^3)\\
    &=(\gamma_j-\lambda_i)\Delta_j-\frac{\gamma_j+\lambda_i}{2}(\gamma_j-\lambda_i)\Delta_j^2+O(\Delta_j^3).
\end{align*}
Since $a_i=O(x)$ and $\Delta_j=O(x)$, substituting these expansions into \eqref{eq:Delta_j} gives
\begin{align*}
    \text{LHS}&=\Delta_j-\frac{a_i}{\lambda_i}(\lambda_i\Delta_j-\frac{1}{2}\lambda_i^2\Delta_j^2+O(\Delta_j^3))=\Delta_j-a_i\Delta_j+\frac{a_i\lambda_i}{2}\Delta_j^2+O(x^3),\\
    \text{RHS}&=\frac{x}{1-r_j}-\frac{1}{1-r_j}\left[r_j\Delta_j-\frac{r_j\gamma_j}{2}\Delta_j^2-r_ja_i\Delta_j+\frac{r_ja_i(\gamma_j+\lambda_i)}{2}\Delta_j^2+O(x^3)\right].
\end{align*}
Equating coefficients of equal powers of $x$ yields the expansion $\Delta_j=x+(\gamma_ir_i+\frac{1}{2}\gamma_j r_j)x^2+O(x^3)$. 

We now return to the cost expressions. Recall that
\begin{align}
    \frac{C_i^{(A)}}{\mu }&=\frac{1}{\mu }(I_i(\eps)+\eps t_i)\nonumber\\
    &=2xt_i-\frac{(1-r_i)t_i^2}{2}-\frac{r_i}{\gamma_i}\lr{t_i-\frac{1-e^{-\gamma_it_i}}{\gamma_i}}.\label{eq:C_iA}
\end{align}
Substituting the expansions of $t_i$ and $1-e^{-\gamma_i t_i}$ into \eqref{eq:C_iA} up to third order (using $t_i^2=x^2+r_i\gamma_ix^3+O(x^4)$, $t_i^3=x^3+O(x^4)$) yields
\begin{align*}
    \frac{C_i^{(A)}}{\mu }=\frac{3}{2}x^2+\frac{2}{3}\gamma_ir_ix^3+O(x^4).
\end{align*}
Similarly, plugging the expansions of $\Delta_j$, $a_i$ and $t_i$ into \eqref{eq:C_iB} and integrating, we obtain
\begin{align*}
    \frac{C_i^{(B)}}{\mu }=\frac{1}{2}x^2+\lr{\frac{1}{2}\gamma_ir_i+\frac{1}{6}\gamma_j r_j}x^3+O(x^4).
\end{align*}
Therefore, the total cost difference between prioritizing class 1 and prioritizing class 2 is
\begin{align*}
    C_1^{(A)}+C_1^{(B)}-C_2^{(A)}-C_2^{(B)}
    &=\mu\left[(\gamma_1 r_1-\gamma_2 r_2)x^3+o(x^3)\right].
\end{align*}
Since $\kappa_1=r_1\gamma_1>r_2\gamma_2=\kappa_2$, this difference is positive as $x\to 0^+$. Hence for small $\eps$ the policy that prioritizes class 2 achieves strictly lower cost than the policy that prioritizes class 1.

\subsubsection{Large Initial Loads \texorpdfstring{$\eps\to \infty$}{eps to inf}}
In this section, we show that when $\eps\to \infty$, the fixed-priority policy prioritizing class 1 yields a strictly lower cost than the policy that always prioritizes class 2. That is, $C_1^{(A)}+C_1^{(B)}<C_2^{(A)}+C_2^{(B)}$ as $\eps\to\infty$.

From \eqref{eq:ti_eps} and the bound $0\leq 1-e^{-y}\leq 1$ for $y\geq 0$,
\begin{align*}
    \frac{\eps}{(1-r_i)\mu }-\frac{r_i}{(1-r_i)\gamma_i}\leq t_i\leq\frac{\eps}{(1-r_i)\mu },
\end{align*}
so as $\eps\to\infty$,
\begin{align*}
    t_i=\frac{\eps}{(1-r_i)\mu }+O(1).
\end{align*}
At time $t_i$, the fraction allocated to class $i$ is
\begin{align*}
    a_i=r_i(1-e^{-\gamma_it_i}),
\end{align*}
with $0\leq r_i-a_i\leq r_ie^{-\gamma_i t_i}=o(1)$. By \eqref{eq:C_iA},
\begin{align*}
    C_i^{(A)}=2\eps t_i-\frac{(1-r_i)\mu}{2} t_i^2-\frac{r_i\mu}{\gamma_i}t_i+\frac{r_i\mu}{\gamma_i^2}(1-e^{-\gamma_i t_i}).
\end{align*}
Since $t_i=\eps/[(1-r_i)\mu]+O(1)$, we can express the stage A cost as
\begin{align*}
    C_i^{(A)}&=\frac{2\eps^2}{(1-r_i)\mu}+O(\eps)-\lr{\frac{\eps^2}{2(1-r_i)\mu}+O(\eps)}-O(\eps)+O(1)\\
    &=\frac{3\eps^2}{2(1-r_i)\mu}+O(\eps).
\end{align*}
Now we continue to derive the Stage B cost. From \eqref{eq:Delta_j} and $q_j^r(\bar{t}_j)=O(1)$ as $\eps\to\infty$, we have
\begin{align}
    \Delta_j=\frac{\eps}{(1-r_j)\mu }+O(1).\label{eq:Delta_j_inf}
\end{align}
Evaluating $q_j^p(t)$ at $t=\bar{t}_j$ and setting it to zero yields
\begin{align*}
    \eps\Delta_j=(1-r_j)\mu \Delta_j^2-\frac{(1-r_j)\mu a_i\Delta_j}{\lambda_i}(1-e^{-\lambda_i\Delta_j})+\frac{r_j\mu \Delta_j}{\gamma_j }(1-e^{-\gamma_j \Delta_j})-\frac{r_j\mu a_i\Delta_j}{\gamma_j -\lambda_i}(e^{-\lambda_i\Delta_j}-e^{-\gamma_j \Delta_j}).
\end{align*}
Substituting this identity for the term $\eps\Delta_j$ in \eqref{eq:C_iB}, and using \eqref{eq:Delta_j_inf}, we derive the Stage B cost
\begin{align*}
    C_i^{(B)}&=\frac{(1-r_j)\mu }{2}\Delta_j^2+\mu \left[\frac{(r_j-1)a_i}{\lambda_i^2}+\frac{r_j}{(\gamma_j )^2}+\frac{r_ja_i}{\lambda_i\gamma_j }\right]+o(1).\\
    &=\frac{1}{2}\cdot\frac{\eps^2}{(1-r_j)\mu }+O(\eps).
\end{align*}
Therefore, as $\eps\to\infty$,
\begin{align*}
    C_i^{(A)}+C_i^{(B)}
    =\frac{\varepsilon^2}{\mu}\left[\frac{3}{2(1-r_i)}+\frac{1}{2(1-r_j)}\right]+o(\varepsilon^2).
\end{align*}
In particular,
\begin{align*}
    \Big(C_2^{(A)}+C_2^{(B)}\Big)-\Big(C_1^{(A)}+C_1^{(B)}\Big)
    =\frac{\varepsilon^2}{\mu}\left(\frac{1}{1-r_2}-\frac{1}{1-r_1}\right)+o(\varepsilon^2).
\end{align*}
Since $r_1<r_2$, the right-hand side is positive for large $\eps$. Hence, as $\eps\to\infty$, the fixed-priority policy that prioritizes class 1 yields a strictly lower cost than the policy that prioritizes class 2.

In conclusion, we showed that when $\kappa_1>\kappa_2$, no fixed-priority policies can be uniformly optimal for all initial states $q_0$.\qed

\subsection{Proofs of Supplementary Lemmas}\label{subsec:proof_supp_lem}
\KiPropsLem*
\textit{Proof. }Recall
\begin{align*}
    K_i(\tau)=\mu \left[(r_i-1)\tau-\frac{r_i}{\gamma_i}(1-e^{-\gamma_i\tau})\right], \tau\geq 0,
\end{align*}
where $r_i\in(0,1)$, $\gamma_i>0$. It follows that $(r_i-1)\tau\leq 0$, $1-e^{-\gamma_i\tau}\geq 0$, $r_i/\gamma_i>0$, so $K_i(\tau)\leq 0$ for all $\tau\geq 0$.

Then, we proceed to prove that $K_1(\tau)\leq K_2(\tau)$ for all $\tau\geq 0$. Define
\begin{align*}
    F(\tau):=\frac{1}{\mu}\left[K_2(\tau)-K_1(\tau)\right]
    =(r_2-r_1)\tau-\frac{r_2}{\gamma_2}(1-e^{-\gamma_2\tau})+\frac{r_1}{\gamma_1}(1-e^{-\gamma_1\tau}).
\end{align*}
Then, $F(0)=0$, and $F'(\tau)=r_2(1-e^{-\gamma_2\tau})-r_1(1-e^{-\gamma_1\tau})$. We aim to show $F'(\tau)\geq 0$ for all $\tau\geq 0$. Observe that $F'(0)=0$, $\lim_{\tau\to\infty}F'(\tau)=r_2-r_1>0$. Differentiating once more yields $F''(\tau)=r_2\gamma_2e^{-\gamma_2 \tau}-r_1\gamma_1e^{-\gamma_1\tau}=\kappa_2 e^{-\gamma_2 \tau}-\kappa_1 e^{-\gamma_1\tau}$. We discuss under the following two cases: 1) $\frac{r_1}{r_2}\cdot\gamma_1\leq \gamma_2\leq \gamma_1$, and 2) $\gamma_2>\gamma_1$.

If $\frac{r_1}{r_2}\cdot\gamma_1\leq \gamma_2\leq \gamma_1$, using $\kappa_2\geq \kappa_1$ and $e^{-\gamma_2\tau}\geq e^{-\gamma_1\tau}$ for $\tau\geq 0$, we have $\kappa_2 e^{-\gamma_2\tau}\geq\kappa_1 e^{-\gamma_2\tau}\geq \kappa_1e^{-\gamma_1\tau}$, so $F''(\tau)\geq 0$ for all $\tau\geq 0$. Therefore, $F'$ is nondecreasing with $F'(0)=0$, hence $F'(\tau)\geq 0$ for all $\tau\geq 0$.

On the other hand, if $\gamma_2>\gamma_1$, we solve $F''(\tau)=0$, obtaining $\tau=s^*:=\log(\kappa_2/\kappa_1)/(\gamma_2-\gamma_1)$. This is nonnegative because $\kappa_2\geq \kappa_1$. Therefore, $F''>0$ on $[0,s^*)$, and $F''<0$ on $(s^*,\infty)$. Hence, $F'$ increases on $[0,s^*]$ and decreases on $[s^*,\infty)$. Together with $F'(0)=0$ and $\lim_{\tau\to\infty}F'(\tau)=r_2-r_1>0$, this implies $F'(\tau)\geq 0$ for all $\tau\geq 0$.

In either case, $F'(\tau)\geq 0$ for all $\tau\geq 0$. Since $F(0)=0$, we obtain $F(\tau)\geq 0$, therefore $K_1(\tau)\leq K_2(\tau)$ for all $\tau\geq 0$.\qed

\nLine \DeltajLem*
\textit{Proof. }Using $0\leq 1-e^{-y}\leq y$ in \eqref{eq:ti_eps}, we have $\mu(1-r_i)t_i\leq \eps$ and $\mu[(1-r_i)t_i+r_it_i]\geq \eps$. Therefore,
\begin{align*}
    x\leq t_i\leq \frac{x}{1-r_i}.
\end{align*}
Then, use $1-e^{-y}\leq y$ on $a_i=\gamma_iq_i^r(t_i)/\mu$, we have
\begin{align*}
    a_i=\frac{\gamma_i}{\mu}\cdot\frac{r_i\mu}{\gamma_i}(1-e^{-\gamma_i t_i})\leq r_i\gamma_i t_i\leq\frac{r_i\gamma_i}{1-r_i}x.
\end{align*}
Therefore, $t_i=O(x)$, $a_i=O(x)$. By \eqref{eq:Delta_j}, we have
\begin{align*}
    \Delta_j-\frac{a_i}{\lambda_i}\lr{1-e^{-\lambda_i\Delta_j}}\leq\frac{\eps}{(1-r_j)\mu}=\frac{x}{1-r_j},
\end{align*}
while
\begin{align*}
    \Delta_j-\frac{a_i}{\lambda_i}(1-e^{-\lambda_i\Delta_j})\geq \Delta_j-a_i\Delta_j=(1-a_i)\Delta_j.
\end{align*}
Combining,  
\begin{align*}
    (1-a_i)\Delta_j&\leq\frac{x}{1-r_j}.
\end{align*}
Because $a_i=O(x)=o(1)$ and $1-a_i=1+o(1)$, $1/(1-a_i)=1+O(a_i)=1+O(x)$. Therefore,
\begin{align*}
    \Delta_j\leq \frac{x}{1-r_j}\cdot\frac{1}{1-a_i}=\frac{x}{1-r_j}(1+O(x))=O(x).
\end{align*}
This proves $\Delta_j=O(x)$.\qed


\end{document}